\documentclass[aop,citesort,MSNbibl,dvips]{arximspdf}
\usepackage{graphics}
\doi{10.1214/08-AOP441}
\volume{37}
\issue{4}
\pubyear{2009}
\firstpage{1502}
\lastpage{1523}

\makeatletter

\newtheorem{lemma}{Lemma}
\newproclaim{Remarks}{Remarks}
\newtheorem{corollary}{Corollary}
\newproclaim{Example}{Example}
\newtheorem{theorem}{Theorem}
\newtheorem{proposition}{Proposition}
\newproclaim{Remark}{Remark}

\makeatother

\begin{document}
\begin{frontmatter}

\title{The structure of the allelic partition of the~total~population
for Galton--Watson processes~with neutral mutations}
\runtitle{Allelic partition of branching processes}

\begin{aug}
\author[A]{\fnms{Jean} \snm{Bertoin}\corref{}\ead[label=e1]{jean.bertoin@upmc.fr}}
\runauthor{J. Bertoin}
\affiliation{Universit\'e Pierre et Marie Curie}
\address[A]{Laboratoire de Probabilit\'es\\
Universit\'e Pierre et Marie Curie\\
175, rue du Chevaleret\\
F-75013 Paris\\
France\\
and\\
DMA\\
Ecole Normale Sup\'erieure\\
Paris\\
France\\
\printead{e1}} 
\end{aug}

\received{\smonth{1} \syear{2008}}
\revised{\smonth{8} \syear{2008}}

%
\begin{abstract}
We consider a (sub-)critical Galton--Watson process with neutral
mutations (infinite alleles model), and decompose the entire population
into clusters of individuals carrying the same allele. We specify the
law of this allelic partition in terms of the distribution of the
number of clone-children and the number of mutant-children of a typical
individual. The approach combines an extension of Harris
representation of Galton--Watson processes and a version of the ballot
theorem. Some limit theorems related to the distribution of the allelic
partition are also given.
\end{abstract}

%
\pdfkeywords{60J80, 60J10,
Branching process,
infinite alleles model,
allelic partition,
ballot theorem}

\end{frontmatter}

\section{Introduction}\label{sec1}

We consider a Galton--Watson process,
that is, a population model with asexual reproduction such that at
every generation, each individual gives birth to a random number of
children according to a fixed distribution and independently of the
other individuals in the population.
We are interested in the situation where a child can be either a clone,
that is,
of the same type (or allele) as its parent, or a mutant, that is, of a
new type.
We stress that each mutant has a~distinct type and in turn gives birth
to clones of itself and to new mutants according to
the same statistical law as its parent, even though it bears a
different allele. In other words, we are working with an infinite
alleles model where mutations are neutral for the population dynamics.
We might as well think of a spatial population model in which children
either occupy the same location as their parents or migrate
to new places and start growing colonies on their own.
This quite basic framework has been often considered in the literature
(see, e.g., \cite{AP0,CD,GP,LyonsPeres,Nerman,Taib});
we also refer to \cite{AD,AP1,AP2,Lambert,LSS,SS} for interesting
variations (these references are of course far from being exhaustive).
Note also that Galton--Watson processes with mutations can be viewed as
a special instance of multitype branching processes (see Chapter V in
Athreya and Ney \cite{AN}
or Chapter 7 in Kimmel and Axelrod \cite{KimmAxel}).

We are interested in the partition of the population into clusters of
individuals having the same allele, which will be referred to as the
\textit{allelic partition}.
Statistics of the allelic partition of a random population model with
neutral mutations
have been first determined in a fundamental work of Ewens \cite{Ewens}
for the Wright--Fisher model
(more precisely this concerns the partition of the population at a
fixed generation).
Kingman \cite{Kingman} provided a deep analysis of this framework, in
connection
with the celebrated coalescent process that depicts the genealogy of
the Wright--Fisher model. We refer to \cite{BG,BBS,DDSJ,DGP,Moehle}
for some
recent developments in this area which involve some related population
models with fixed generational size and certain exchangeable coalescents.

The main purpose of the present work is to describe explicitly the
structure of the allelic partition of the entire population
for Galton--Watson processes with neutral mutations.
We will always assume that the Galton--Watson process is critical or
subcritical, so the descent of any individual becomes eventually
extinct, and in particular the allelic clusters are finite a.s.
We suppose that every ancestor (i.e., individual in the initial
population) bears a different allele; it is convenient to
view each ancestor as a mutant of the zeroth kind. We then call mutant
of the first kind a mutant-child of
an individual of the allelic cluster of an ancestor, and the set of all
its clones (including that mutant) a cluster of the first kind. By iteration,
we define mutants and clusters of the $k$th kind for any integer
$k\geq0$.

In order to describe the statistics of the allelic partition,
we distinguish an ancestor which will then be referred to as \textit{Eve},
and focus on its descent.
The set of all individuals bearing the same allele as Eve is called the
\textit{Eve cluster}.
The Eve cluster has obviously the genealogical structure of a
Galton--Watson tree
with reproduction law given by the distribution of the number of
clone-children of a typical individual.
Informally, the branching property indicates that the same holds for
the other clusters of the allelic partition. Further,
it should be intuitively clear that the process which counts the number
of clusters of the $k$th kind for $k\geq0$ is again
a~Galton--Watson process whose reproduction law is given by the
distribution of the number of mutants of the first kind; this
phenomenon has already been pointed at in the work of Ta\"{\i}b \cite{Taib}.
That is to say that, in some loose sense the allelic partition
inherits branching structures from the initial Galton--Watson process.
Of course, these formulations are only heuristic and precise statements
will be given later on. We also stress that
the forest structure which connects clusters of different kinds and the
genealogical structure on each cluster are not independent since, typically,
the number of mutants of the first kind who stem from the Eve cluster
is statistically related to the size of the Eve cluster.

Our approach essentially relies on a variation of the well-known
connection due to Harris \cite{Harris1,Harris2} between ordinary
Galton--Watson processes and sequences of i.i.d. integer-valued random
variables. Specifically, we incorporate neutral mutations in Harris
representation and by combination with the celebrated ballot theorem
(which is another classical tool in this area as it is expounded,
e.g., by Pitman; see Chapter 6 in \cite{PiSF}), we obtain expressions
for the joint distribution of various natural variables (size of the
total descent of an ancestor, number of alleles, size and number of
mutant-children of an allelic cluster) in terms of the transition
probabilities of the two-dimensional random walk which is generated by
the numbers of clone-children and of mutant-children of a typical individual.

We also investigate some limit theorems in law; typically we show that
when the numbers
of clone-children and mutant-children of an individual are independent
(and some further technical conditions), the sequence of the relative
sizes of the allelic clusters in a typical tree has a limiting
conditional distribution when the size
of the tree and the number of types both tend to infinity according to
some appropriate regime. The limiting distribution that arises has
already appeared in the study of the standard additive coalescent by
Aldous and Pitman \cite{AP1}.
We also point at limit theorems for allelic partitions of Galton--Watson
forests, where, following Duquesne and Le Gall \cite{DuLG,DuLG2}, the
limits are described in terms of certain L\'evy trees. In particular,
this provides an explanation to a rather striking identity between two
self-similar fragmentation processes that were defined on the one hand
by logging the Continuum Random Tree
according to a Poisson point process along its skeleton \cite{AP1}, and
on the other hand by splitting the unit-interval at instants when the
standard Brownian excursion with a negative drift reaches new infima
\cite{Ber1}.

\section{Allelic partitions in a Galton--Watson forest}\label{sec2}

We first develop some material and notation about Galton--Watson
forests with neutral mutations, referring to Chapter 6 in Pitman \cite
{PiSF} for background in the case without mutations.

\subsection{Basic setting}\label{sec21}
Let
\[
\xi=\bigl(\xi^{(\mathrm{c})},\xi^{(\mathrm{m})}\bigr)
\]
be a pair of nonnegative integer-valued random variables which should
be thought of respectively as the number of clone-children and the
number of mutant-children of a typical individual. We also write
\[
\xi^{(+)}=\xi^{(\mathrm{c})}+\xi^{(\mathrm{m})}
\]
for the total number of children, and assume throughout this work that
\[
\mathbb{E}\bigl(\xi^{(+)}\bigr)\leq1 ,
\]
that is, we work in the critical or subcritical regime.
We implicitly exclude the degenerate case when $\xi^{(\mathrm{c})}\equiv0$
or $\xi
^{(\mathrm{m})}\equiv0$
and, as a consequence, the means $\mathbb{E}(\xi^{(\mathrm{c})})$ and
$\mathbb{E}(\xi^{(\mathrm{m})}
)$ are
always less than $1$.

We write $\mathbb{Z}_+$ and $\mathbb{N}$ for the sets of nonnegative
integers and
positive integers, respectively.
A pair $(g,n)\in\mathbb{Z}_+\times\mathbb{N}$ is then used to
identify an individual
in an infinite population model, where the first coordinate $g$ refers
to the generation and the second coordinate $n$ to the rank of the
individual of that generation
(we stress that each generation consists of an infinite sequence of
individuals).
We assume that each individual at generation $g+1$ has a unique parent
at generation $g$.
We consider a family
\[
(\xi_{g,n}\dvtx g\in\mathbb{Z}_+\mbox{ and } n\in\mathbb{N})
\]
of i.i.d. copies of $\xi$ which we use to define the Galton--Watson
process with neutral mutations. Specifically, $\xi_{g,n}=(\xi^{(\mathrm{c})}
_{g,n},\xi^{(\mathrm{m})}_{g,n})$ is the pair given by the number of clone-children
and mutant-children of the $n$th individual at generation $g$. We may
assume that the offspring of each individual is ranked, which induces a
natural order at the next generation by requiring further that if
$(g,n)$ and $(g,n')$ are two individuals at the same generation $g$
with $n<n'$, then at generation $g+1$ the children of $(g,n)$ are all
listed before those of $(g,n')$.

\subsection{Encoding the Galton--Watson forest with mutations}\label{sec22}

Next, we enumerate as follows the individuals of the entire population
(i.e., of all generations) by a~variation of the well-known depth-first
search algorithm that takes mutations into account.
We associate to each individual a label $(a,m,s)$, where $a\in\mathbb
{N}$ is
the rank of the ancestor in the initial population, $m$ the number of
mutations and $s$ a finite sequence of positive integers which keeps
track of the genealogy of the individual. Specifically, the label of
the $a$th individual in the initial generation $g=0$ is
$(a,0,\varnothing)$. If an individual at the $g$th generation has the
label $(a,m,(i_1,\ldots, i_g))$,
and if this individual has $j^{(\mathrm{c})}$ clone-children and $j^{(\mathrm{m})}$
mutant-children, then the labels assigned to its clone-children are
\[
(a,m, (i_1,\ldots, i_g, 1)), \ldots, \bigl(a,m, \bigl(i_1,\ldots, i_g, j^{(\mathrm{c})}\bigr)\bigr) ,
\]
whereas the labels assigned to its mutant-children are
\[
\bigl(a,m+1, \bigl(i_1,\ldots, i_g, j^{(\mathrm{c})}+1\bigr)\bigr), \ldots, \bigl(a,m+1,
\bigl(i_1,\ldots, i_g,j^{(\mathrm{c})}
+ j^{(\mathrm{m})}\bigr)\bigr) .
\]

\begin{figure}

\includegraphics{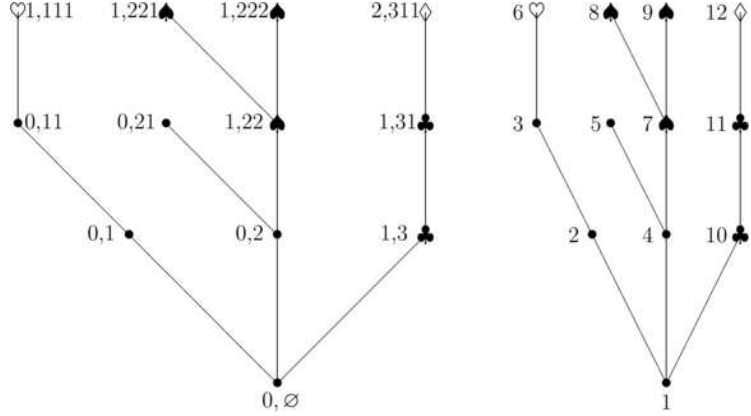}

\caption{Depth-first search with mutations on a
genealogical tree.
The symbols $\bullet, \spadesuit, \heartsuit, \Diamond, \clubsuit$
represent the different alleles.
Left: the label $(m,s)$ of an individual is given by the number $m$ of
mutations and the sequence $s$ that specifies its genealogy; for the
sake of simplicity, the rank $a$ of the ancestor has been omitted.
Right: the same tree with individuals ranked by the depth-first search
algorithm with mutations.}\label{fig1}
\end{figure}

Clearly, any two distinct individuals have different labels.
We then introduce the (random) map
\[
\rho\dvtx\mathbb{N}\to\mathbb{Z}_+\times\mathbb{N},
\]
which consists in ranking the individuals in the lexicographic order of
their labels;
see Figure \ref{fig1}.
That is to say that $\rho(i)=(g,n)$ if and only if the $i$th
individual in the lexicographic order of labels corresponds to the
$n$th individual at generation~$g$.
This procedure for enumerating the individuals will be referred to as
the \textit{depth-first search algorithm with mutations}. We shall
also use the notation
\[
\xi_i=\xi_{\rho(i)} ,\qquad i\in\mathbb{N} ,
\]
and whenever no generation is specified, the terminology $i$th
individual will implicitly refer to the rank of that individual induced
by depth-first search with mutation, that is, the $i$th individual
means the $n$th individual at generation $g$ where $\rho(i)=(g,n)$.

%

\begin{lemma} \label{Le1}
\textup{(i)} The variables $\xi_1, \xi_2,
\ldots
$ are i.i.d. with the same law as $\xi$.

\textup{(ii)} The sequence $(\xi_{g,n}\dvtx g\in\mathbb{Z}_+
\mbox{ and }
n\in\mathbb{N})$ can be recovered from $(\xi_i\dvtx i\in\mathbb{N})$ a.s.
\end{lemma}

\begin{pf} It should be plain from the definition of the depth-first
search algorithm with mutations that for every $i\in\mathbb{N}$,
$\rho(i+1)$ is
a deterministic function of $\xi_1, \ldots, \xi_i$ which takes
values in
$(\mathbb{Z}_+\times\mathbb{N})\setminus\{\rho(1), \ldots, \rho
(i)\}$. Since
$(\xi
_{g,n}\dvtx g\in\mathbb{Z}_+ \mbox{ and } n\in\mathbb{N})$ is a
sequence of i.i.d.
variables with the same law as $\xi$, this yields the first claim by
induction. The second claim follows from the fact that each individual
has a finite descent a.s. [because the Galton--Watson process is \mbox{(sub-)}critical],
which easily entails that the map $\rho$ is bijective.
Further, it is readily seen that the inverse bijection
is a function of the sequence $(\xi_i\dvtx i\in\mathbb{N})$.
\end{pf}

Henceforth, we shall therefore encode the Galton--Watson process with
neutral mutations by a sequence $(\xi_i\dvtx i\in\mathbb{N})$ of i.i.d.
copies of
$\xi$. We denote by
$({\mathcal F}_i)_{i\in\mathbb{N}}$ the natural filtration generated
by this sequence.

We next briefly describe the genealogy of the Galton--Watson process as
a forest of i.i.d. genealogical trees. Denote for every $n\in\mathbb
{N}$ by
\[
\alpha_n=\rho^{-1}(0,n),
\]
so that $\alpha_1=1 < \alpha_2< \cdots$ is the increasing sequence of
the ranks of ancestors induced by the depth-first search algorithm with
mutations.
For example, \mbox{$\alpha_2=13$} in the situation described by Figure \ref
{fig1}. The
procedure for labeling individuals ensures that the descent of the
$i$th ancestor $\alpha_i$ corresponds
to the integer interval
\[
{[\alpha_i, \alpha_{i+1}[}:=\{\alpha_i, \alpha_i+1, \ldots, \alpha
_{i+1}-1\}
\]
(that is to say, if we index the population model using generations,
then the descent of
$(0,i)$ is the image of $[\alpha_i, \alpha_{i+1}[$ by the inverse
bijection $\rho^{-1}$).

We write
\[
\mathbb{T}_i:=(\xi_{\alpha_i-1+\ell}\dvtx1\leq\ell\leq\alpha
_{i+1}-\alpha_i)
\]
for the finite sequence of the numbers of clone-children and
mutant-children of the individuals in the descent of the $i$th
ancestor. So $\mathbb{T}_i$ encodes
(by the depth-first search algorithm with mutations) the genealogical
tree of the $i$th ancestor, and it should be intuitively clear that
the family
$(\mathbb{T}_i\dvtx i\in\mathbb{N})$ is a forest consisting in a
sequence of i.i.d.
genealogical trees. To give a rigorous statement, it is convenient to
introduce the downward skip-free (or left-continuous) random walk
%
%
\begin{equation}\label{ES+}
S^{(+)}_n:=\xi^{(+)}_1+\cdots+\xi^{(+)}_n -n ,\qquad n\in\mathbb{Z}_+ ,
\end{equation}
and the passage times
%
%
\begin{equation}\label{ET+}
T^{(+)}_i:=\inf\bigl\{n\geq0\dvtx S^{(+)}_n=-i\bigr\} ,\qquad i\in\mathbb{Z}_+ .
\end{equation}
We stress that the $T^{(+)}_i$ form an increasing sequence of
$({\mathcal F}
_n)$-stopping times.

\begin{lemma} \label{Le2} There is identity
\[
\alpha_i- 1= T^{(+)}_{i-1}
\]
for every $i\in\mathbb{N}$ and, as a consequence, the sequence
$\mathbb{T}_1, \ldots$
is i.i.d.
\end{lemma}

\begin{pf} This formula is a close relative of the classical identity
of Dwass \cite{Dwass} and would be well known if individuals were
enumerated by the usual depth-first search algorithm (i.e., without
taking care of mutations), see, for example, Lemma~63 in \cite{PiSF} or
\cite{LG}. The proof in the present case is similar. Indeed the formula
is obvious for $i=1$, and for $i=2$,
we have on the one hand that
\[
\alpha_2-1=1+\xi^{(+)}_1
+\cdots+\xi^{(+)}_{\alpha_2-1}
\]
by expressing the fact that
the predecessor of the second ancestor found by depth-first search with
mutations has a rank given by the size of the population generated by
Eve, that is,
Eve herself and her descendants. On the other hand,
we must have $1+\xi^{(+)}_1
+\cdots+\xi^{(+)}_n>n$ when $n<\alpha_2-1$, since otherwise the
depth-first search algorithm with mutations would explore the second
ancestor before having completed the exploration of the entire descent
of Eve. This proves the identity for $i=2$, and the general case then
follows by iteration. Finally, the last claim is an immediate
consequence of Lemma \ref{Le1}(i) and the strong Markov property.
\end{pf}

\subsection{Allelic partitions}\label{sec23}

We can now turn our attention to defining allelic partitions.
In this direction, recall that every ancestor has a different type
(i.e., bears a different allele), and thus should be viewed as an
initial mutant. More generally, we call \textit{mutant}
an individual which either belongs to the initial generation or is the
mutant-child of some individual, and then write
\[
1=\mu_1<\mu_2<\cdots
\]
for the ranks of mutants in the depth-first search algorithm with mutations.
For example, $\mu_2=6$, $\mu_3=7$, $\mu_4=10$, $\mu_5=12$ and $\mu
_6=\alpha_2=13$ in the situation depicted by Figure \ref{fig1}.
The upshot of this algorithm is that the set of individuals that bear
the same allele as the $j$th mutant $\mu_j$ corresponds precisely to
the integer interval
$[\mu_j, \mu_{j+1}[$. In this direction, it is therefore natural to
introduce for every
$j\in\mathbb{N}$ the
$j$th \textit{allelic cluster}
\[
\mathbb{C}_j:=( \xi_{\mu_j-1+\ell}\dvtx1\leq\ell\leq\mu_{j+1}-\mu
_j) ,
\]
that is, $\mathbb{C}_j$ is the finite sequence of the numbers of clone-children
and mutant-children of the individuals bearing the same allele as the
$j$th mutant. The sequence $(\mathbb{C}_j)_{j\in\mathbb{N}}$
encodes the allelic
partition of the entire population.

\begin{Remarks*}
1. Each allelic cluster $\mathbb{C}_j$ is naturally endowed with a
structure of
rooted planar tree which is induced by the Galton--Watson process. More
precisely, the latter is encoded via the usual depth-first search
algorithm by the sequence $(\xi^{(\mathrm{c})}_{\mu_j-1+\ell}\dvtx1\leq\ell
\leq\mu
_{j+1}-\mu_j)$; in particular the $j$th mutant $\mu_j$ is viewed as
the root
(i.e., ancestor) of the cluster $\mathbb{C}_j$. In other words, the depth-first
search algorithm with mutations
for the Galton--Watson process
induces precisely the usual depth-first search applied to the forest of
allelic clusters viewed as a sequence of planar rooted trees.

2. We also stress that the initial Galton--Watson process can
be recovered from the allelic partition $(\mathbb{C}_j)_{j\in\mathbb
{N}}$. Indeed, the
previous observation shows how to construct the portion of the
genealogical tree corresponding to the allelic cluster generated by an
initial mutant, and the latter also contains the information which is
needed to identify the mutant-children of the first kind.
Mutant-children of the first kind are the roots of the subtrees
corresponding to the allelic clusters of the second kind, and by
iteration the entire genealogical forest can be recovered.

Just as above, it is now convenient to introduce the downward skip-free
random walk
%
%
\begin{equation}\label{ESc}
S^{(\mathrm{c})}_n:=\xi^{(\mathrm{c})}_1+\cdots+\xi^{(\mathrm{c})}_n -n ,\qquad
n\in\mathbb{Z}_+ ,
\end{equation}
and the passage times
%
%
\begin{equation}\label{ETc}
T^{(\mathrm{c})}_j:=\inf\bigl\{n\geq0\dvtx S^{(\mathrm{c})}_n=-j\bigr\} ,\qquad j\in\mathbb
{Z}_+ .
\end{equation}
Again, the $T^{(\mathrm{c})}_j$ form an increasing sequence of $({\mathcal F}
_i)$-stopping times.
\end{Remarks*}

\begin{lemma}\label{Le3} There is identity
\[
\mu_j- 1=T^{(\mathrm{c})}_{j-1}
\]
for every $j\in\mathbb{N}$. As a consequence, for every $j\in\mathbb
{N}$, $\mathbb{C}_j$ is
adapted to the sigma-field
${\mathcal F}_{T^{(\mathrm{c})}_{j}}$, whereas $\mathbb{C}_{j+1}$ is
independent of ${\mathcal F}
_{T^{(\mathrm{c})}_{j}}$ and
has the same distribution as~$\mathbb{C}_1$. In particular
the sequence of the allelic clusters $\mathbb{C}_1, \mathbb{C}_2,
\ldots$ is i.i.d.
\end{lemma}

The proof is similar to that of Lemma \ref{Le2} and therefore omitted.

We also introduce the number of alleles, that is, of different types,
which are present in the $i$th tree $\mathbb{T}_i$:
\[
{A}_i:=\operatorname{Card}\{j\in\mathbb{N}\dvtx\mu_j\in[\alpha_i,\alpha
_{i+1}[\} ;
\]
for example, $A_1=5$ in the situation described by Figure \ref{fig1}.
Note that there is the alternative expression
\[
{A}_i=1+\sum_{\alpha_i\leq\ell< \alpha_{i+1}}\xi^{(\mathrm{m})}_{\ell}
.
\]

\begin{corollary} \label{C1}
\textup{(i)} For every $i\in\mathbb{Z}_+$,
we have
\[
\alpha_{i+1}=\mu_{{A}_1+\cdots+ {A}_i+1} ;
\]
equivalently, there is the identity
\[
T^{(+)}_i=T^{(\mathrm{c})}_{{A}_1+\cdots+{A}_i} .
\]

\textup{(ii)} The allelic partition of the tree $\mathbb{T}_i$,
which is
induced by
restricting the allelic partition of the entire population to $\mathbb{T}_i$,
is given by
\[
(\mathbb{C}_{{A}_1+\cdots+{A}_{i-1}+\ell}\dvtx1\leq\ell\leq{A}_i) .
\]
As a consequence, the sequence of the allelic partitions of the trees
$\mathbb{T}_i$ for $i\in\mathbb{N}$,
is i.i.d.
\end{corollary}

\begin{pf} (i) The first identity should be obvious from the definition
of the depth-first search with mutations, as ${A}_1+\cdots+ {A}_i$ is
the number of alleles which have been found after completing the
exploration of the $i$ first trees and the next mutant is then the
$(i+1)$th ancestor. The second then follows from Lemmas~\ref{Le2}~and~\ref{Le3}.

(ii) The first assertion is immediately seen from (i) and the
definitions of the trees and of the allelic clusters. Then observe that
the number ${A}_i$ of alleles in the tree $\mathbb{T}_i$
is a function of that tree, and so is the allelic partition. The second
assertion thus derives from Lemma \ref{Le2}.
\end{pf}

It may be interesting to point out that $(T^{(+)}_{i}, i\geq0)$ and
$(T^{(\mathrm{c})}
_j, j\geq0)$ are both increasing random walks. The range
\[
{\mathcal R}^{(+)}:= \bigl\{T^{(+)}_i\dvtx i\geq0 \bigr\}
\]
is the set of predecessors of ancestors (in the depth-first search
algorithm with mutations),
whereas
\[
{\mathcal R}^{(\mathrm{c})}:= \bigl\{T^{(\mathrm{c})}_j\dvtx j\geq0 \bigr\}
\]
corresponds to predecessors of mutants. These
are two regenerative subsets of $\mathbb{Z}_+$, in the sense that each
can be viewed
as the set of renewal epochs of some recurrent event (cf. Feller \cite
{Fel1,Fel2}).
Observe that both yield a partition of the set of positive integers
into disjoint intervals:
\[
\mathbb{N}=\bigcup_{i\geq1} \bigl]T^{(+)}_{i-1},T^{(+)}_{i}\bigr] =\bigcup_{j\geq
1} \bigl]T^{(\mathrm{c})}
_{j-1},T^{(\mathrm{c})}_{j}\bigr],
\]
that correspond respectively to the trees in the Galton--Watson forest
and to the allelic clusters.
By Corollary \ref{C1}(i), there is the embedding
\[
{\mathcal R}^{(+)}\subseteq{\mathcal R}^{(\mathrm{c})}
\]
and more precisely, this embedding is compatible with regeneration, in
the sense
that for every $k\in\mathbb{Z}_+$, conditionally on $k\in{\mathcal
R}^{(+)}$,
the shifted sets
${\mathcal R}^{(+)}\circ\theta_k:=\{i\geq0\dvtx k+i\in{\mathcal
R}^{(+)}\}$ and
${\mathcal R}^{(\mathrm{c})}\circ\theta
_k:=\{j\geq0\dvtx k+j\in{\mathcal R}^{(\mathrm{c})}\}$
are independent of the sigma-field ${\mathcal F}_k$ generated by $(\xi_1,
\ldots,
\xi_k)$ and their joint law is the same as that
of $({\mathcal R}^{(+)},{\mathcal R}^{(\mathrm{c})})$. We refer to \cite{Ber1}
for applications of this notion.
Roughly speaking, this implies that the allelic split of each interval
$]T^{(+)}_{i-1},T^{(+)}_{i}]$
produces smaller intervals $]T^{(\mathrm{c})}_{j-1},T^{(\mathrm{c})}_{j}]$ in a
random way that only
depends on the length $T^{(+)}_{i}-T^{(+)}_{i-1}$ (i.e., the size of
$\mathbb{T}
_{i}$), independently of its location and of the other integer
intervals. This can be thought of as a fragmentation property (see
\cite
{RFCP}) for the sizes of the trees.

\subsection{Allelic trees and forest}\label{sec24}

In order to analyze the structure of allelic partitions, we introduce
some related notions.
The genealogy of the population model naturally induces a structure of
forest on the set of different alleles. More precisely, we enumerate
this set by declaring that the $j$th allele is that of the $j$th
cluster $\mathbb{C}_j$, and define a planar graph on the set of
alleles (which
is thus identified as $\mathbb{N}$) by drawing an edge between two
integers $j
< k$
if and only if the parent of the $k$th mutant $\mu_k$ is an individual
of the $j$th allelic cluster $\mathbb{C}_j$. This graph is clearly a forest
(i.e., it contains no cycles), which we call
the \textit{allelic forest}, and more precisely the $i$th allelic tree is
that induced by the
mutant descent of the $i$th ancestor $\alpha_i$. In other words, the
$i$th allelic tree
is the genealogical tree of the different alleles present in $\mathbb
{T}_i$. In
particular,
the sequence of allelic trees is i.i.d. and their sizes are given by
$({A}_i, i\in\mathbb{N})$.

Recall that the \textit{breadth-first search} in a forest consists in
enumerating individuals in the lexicographic order of their labels,
where the label of the $n$th individual at generation $g$ is now given
by the triplet $(a,g,n)$, with $a$ the rank of the ancestor at the
initial generation. After a (short) moment of thought, we see that the
definition of depth-first search with mutations for the Galton--Watson
process ensures that the labeling of alleles by integers agrees with
breadth-first search on the allelic forest, in the sense that the
$j$th allele is found at the $j$th step of the breadth-first search
on the allelic forest.

\begin{figure}[b]

\includegraphics{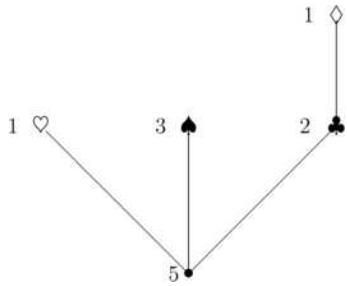}

\caption{Allelic tree corresponding to the
genealogical tree with mutations in Figure \protect\ref{fig1}. The
labels represent the
sizes of the allelic clusters.}
\label{fig2}
\end{figure}

For every $j\in\mathbb{N}$, we consider the number of new mutants who are
generated by the $j$th allelic cluster, viz.
\[
M_j:=\sum_{\mu_j\leq\ell< \mu_{j+1}}\xi^{(\mathrm{m})}_{\ell} .
\]
For instance, we have $M_1= 3$, $M_4=1$ and $M_2=M_3=M_5=0$ in the
situation depicted by Figures \ref{fig1} and \ref{fig2}.
The allelic forest is thus encoded by breadth-first search via the
sequence $ (M_{j}, j\in\mathbb{N} )$.

%

\begin{lemma} \label{Le4} The sequence $(M_j, j\in\mathbb{N})$ is
i.i.d., and
therefore the
allelic forest is a Galton--Watson forest with reproduction law the
distribution of $M_1$.
As a~consequence, the size ${A}_1$ of the first allelic tree
is given by the identity
\[
{A}_1=\min\{j\geq1\dvtx M_1+\cdots+M_j=j-1\} ,
\]
showing that ${A}_1$ is an $({\mathcal F}_{T^{(\mathrm{c})}_{j}})$-stopping time.
\end{lemma}

\begin{pf} Recall from Lemma \ref{Le3} that the sequence $\mathbb
{C}_1, \mathbb{C}_2,
\ldots$ of the allelic clusters is i.i.d. Clearly, each variable $M_j$
only depends on $\mathbb{C}_j$, which entails our first claim. The second
follows from the well-known fact that breadth-first search induces
a bijective transformation between the distributions of (sub-)critical
Galton--Watson forests and those
of i.i.d. sequences of integer-valued variables with mean less than or
equal to one
(see, e.g., Section 6.2 in \cite{PiSF}).

Finally, the identity for the number ${A}_1$ of alleles present in the
tree $\mathbb{T}_1$
follows from the preceding observations and again a~variation of the
celebrated formula of Dwass \cite{Dwass} (see, e.g., Lemma \ref{Le2} in
the present work), as plainly, ${A}_1$ coincides with the total size of
the first tree in the allelic forest.
\end{pf}

\section{Some applications of the ballot theorem}\label{sec3}

We start by stating a version of the classical ballot theorem that
will be used in this section; see \cite{Takacs}.
Let $(X_1,\ldots, X_{n})$ be an $n$-tuple of random variables with
values in some space
$E$, which is cyclically exchangeable, in the sense that for every
$i\in\mathbb{N}$,
there is the identity in law
\[
(X_1,\ldots, X_{n})
\stackrel{\mathcal{L}}{=}
(X_{i+1},\ldots, X_{i+n}) ,
\]
where we agree that addition of indices is taken modulo $n$.
Consider a function
\[
f\dvtx E\to\{-1, 0, 1, 2, \ldots\}
\]
and assume that
\[
\sum_{j=1}^n f(X_j)=-k
\]
for some $1\leq k \leq n$.

\begin{lemma}[(Ballot theorem)]\label{Le6}
Under the assumptions above, the probability that the process of the
partial sums of the sequence $f(X_1), \ldots, f(X_n)$ remains above
$-k$ until the $n$-step
is
\[
\mathbb{P} \Biggl( \min\Biggl\{j\geq1\dvtx\sum_{i=1}^j f(X_i)=-k\Biggr\}=n
\Biggr)=k/n .
\]
\end{lemma}

\subsection{Distribution of the allelic tree}\label{sec31}

We have now introduced all the tools which are needed for describing
some statistics of the allelic partition of a Galton--Watson tree with
neutral mutations. We only need one more notation. We write
%
%
\begin{equation}\label{E9}\pi_{k,\ell}=\mathbb{P}\bigl(\xi^{(\mathrm{c})}=k,\xi^{(\mathrm{m})}
=\ell\bigr) ,\qquad
k,\ell\in\mathbb{Z}_+,
\end{equation}
for the probability function of the reproduction law of the
Galton--Watson process with mutations.
For every integer $n\geq1$, we also write $\pi^{*n}$ for the $n$th
convolution product of that law, that is,
\[
\pi^{*n}_{k,\ell}=\mathbb{P}\bigl(\xi^{(\mathrm{c})}_1+\cdots+\xi^{(\mathrm{c})}_n=k,
\xi^{(\mathrm{m})}_1+\cdots+\xi^{(\mathrm{m})}
_n=\ell\bigr).
\]

\begin{Example*}
Suppose that the dynamics of the population
can be described as follows. We start from a usual Galton--Watson
process with reproduction law on~$\mathbb{Z}_+$, say $\varrho$, and
assume that
at each step mutations affect each child with probability $p\in\,]0,1[$,
independently of the other children.
In other words, the allelic forest is obtained by pruning or
percolation on the genealogical forest of the Galton--Watson process,
cutting each edge with probability $p$ and independently of the other
edges. See, for example, Aldous and Pitman \cite{AP0} or Chapter 4 in
Lyons and Peres \cite{LyonsPeres}.
Analytically, this means that if $\xi$ is a random variable with law
$\varrho$, then the conditional distribution of $(\xi^{(\mathrm{c})},\xi
^{(\mathrm{m})})$ given
$\xi=k$ is that of $(k-B(k,p), B(k,p))$,
where $B(k,p)$ denotes a binomial variable with parameters $k$ and $p$.
In this situation, it is easily seen that
%
%
\begin{equation}\label{E3}
\pi^{*n}_{k,\ell}= \pmatrix{ k+\ell\cr k}
(1-p)^k p^\ell\varrho^{*n}_{k+\ell}
\end{equation}
with $\varrho^{*n}$ denoting the $n$th convolution power of $\varrho$.
This expression is entirely explicit when $\varrho$ is, for example,
the Poisson, or binomial or geometric, distribution as in those cases,
there are known formulas for $\varrho^{*n}$.
Of course, there are other natural examples in which the
two-dimensional probability function $\pi^{*n}$
can be expressed in terms of simpler one-dimensional probability
functions, for instance,
when $\xi^{(\mathrm{c})}$ and $\xi^{(\mathrm{m})}$ are assumed to be independent or
when $\xi^{(\mathrm{c})}
=\beta\xi$
and $\xi^{(\mathrm{m})}=(1-\beta)\xi$ where $\beta$ stands for a Bernoulli variable
which is independent of $\xi$.
\end{Example*}

Corollary \ref{C1} enables us to restrict our attention to the allelic
partition of the tree generated by a typical ancestor, say for
simplicity, Eve.
Recall that $T^{(+)}_1$ denotes the size of the genealogical tree
$\mathbb
{T}_1$ of
Eve, that ${A}_1$ is the number of alleles found in $\mathbb{T}_1$
and that the $j$th allelic cluster $\mathbb{C}_j$ generates $M_j$
mutant-children. Further, we know from Lemma \ref{Le4} that the first
allelic tree is encoded by breadth-first search via the finite sequence
$ (M_{j}, 1\leq j \leq{A}_1 )$. The latter only retains
partial information about the structure of the allelic partition of
$\mathbb{T}
_1$, and thus it is natural to enrich it by considering more generally
the sequence of pairs $ ((|\mathbb{C}_j|,M_{j}), 1\leq j \leq{A}_1
)$, where
\[
|\mathbb{C}_j|:= \mu_{j+1}-\mu_j
\]
denotes the size of the $j$th allelic cluster, that is, the number of
individuals having the $j$th type. In other words, we enrich the
allelic tree
by assigning to each allele the size of the corresponding allelic cluster.
We may now state our main result, which can be viewed as a
generalization of a celebrated identity due to Dwass~\cite{Dwass}.

\begin{theorem} \label{T1}
\textup{(i)} The joint law of the size of $\mathbb{T}_1$ and
its number
of alleles is given by
\[
\mathbb{P}\bigl(T^{(+)}_1=n, {A}_1=k\bigr)=\frac{1}{n}\pi^{*n}_{n-k,k-1}
,\qquad
1\leq k
\leq n .
\]

\begin{longlist}[(iii)]
\item[(ii)] The joint law of the size of the Eve cluster and
the number of its mutant-children is given by
\[
\mathbb{P}(|\mathbb{C}_1|=n, M_1=\ell)=\frac{1}{n}\pi
^{*n}_{n-1,\ell} ,\qquad
n\geq
1\mbox{ and } \ell\geq0 .
\]

\item[(iii)]
For every integers $k\geq1$, $n_1, \ldots, n_k \geq1$ and $\ell_1,
\ldots, \ell_k\geq0$ such that
\[
\sum_{i=1}^j \ell_i>j-1\qquad \mbox{whenever } 1\leq j < k ,
\]
we have
\[
\mathbb{P}(|\mathbb{C}_1|=n_1, M_1=\ell_1, \ldots
, |\mathbb{C}_k|=n_k, M_k=\ell_k)
= \prod_{i=1}^k\frac{1}{n_i}\pi^{*n_i}_{n_i-1,\ell_i} .
\]
\end{longlist}
\end{theorem}

\begin{Remarks*}
1. Restricting our attention in part (iii) to sequences\break $\ell_1,
\ldots
, \ell_k\geq0$ with
\[
\inf\Biggl\{j\geq1\dvtx\sum_{i=1}^j \ell_i=j-1\Biggr\}= k ,
\]
we stress that the statement describes the law of the entire allelic tree.

2. In particular, the law of the number ${A}_1$ of alleles is given by
\[
\mathbb{P} ({A}_1=k ) = \sum_{n=1}^{\infty}n^{-1}\pi
^{*n}_{n-k,k-1}
,\qquad k\geq0 .
\]
It may be interesting to point out that there is also the formula
\[
\mathbb{P} ({A}_1=k ) =\frac{1}{k}\nu^{*k}_{k-1},
\]
where
\[
\nu_{\ell}=\mathbb{P}(M_1=\ell)=\sum_{n=1}^{\infty} n^{-1}\pi
^{*n}_{n-1,\ell}
\]
and $\nu^{*k}$ the $k$th convolution power of $\nu$. Indeed, this alternative
formulation is seen from Lemma \ref{Le4} and Dwass formula \cite{Dwass}.
\end{Remarks*}

\begin{pf} Recall that $(\xi^{(\mathrm{c})}_1, \xi^{(\mathrm{m})}_1), \ldots, (\xi
^{(\mathrm{c})}_n, \xi^{(\mathrm{m})}_n)$
is a sequence
of $n$ i.i.d. copies of $(\xi^{(\mathrm{c})},\xi^{(\mathrm{m})})$ and consider the
partial sums
of coordinates
\[
\Sigma^{(\mathrm{c})}_j=\sum_{i=1}^j\xi^{(\mathrm{c})}_i ,\qquad \Sigma
^{(\mathrm{m})}_j=\sum_{i=1}^j\xi^{(\mathrm{m})}_i
\quad\mbox{and}\quad \Sigma_j=\Sigma^{(\mathrm{c})}_j+\Sigma^{(\mathrm{m})}_j .
\]
Introduce for every $1\leq k \leq n$ the event
\[
\Lambda_{n-k,k-1}= \bigl\{\Sigma^{(\mathrm{c})}_n=n-k,\Sigma^{(\mathrm{m})}
_n=k-1 \bigr\}= \bigl\{
\Sigma_n=n-1, \Sigma^{(\mathrm{m})}_n=k-1 \bigr\}
\]
and observe that the sequence
$(\xi^{(\mathrm{c})}_1, \xi^{(\mathrm{m})}_1), \ldots, (\xi^{(\mathrm{c})}_n, \xi^{(\mathrm{m})}_n)$
is (cyclically)
exchangeable conditionally on
$\Lambda_{n-k,k-1}$. Further, we have by definition that
\[
\mathbb{P}(\Lambda_{n-k,k-1})=\pi^{*n}_{n-k,k-1} .
\]

Plainly, there is the identity
\[
\bigl\{T^{(+)}_1=n, {A}_1=k\bigr\}=\Lambda_{n-k,k-1}\cap
\bigl\{\min\{j\geq1\dvtx\Sigma_j=j-1\}=n\bigr\}
\]
as, according to Lemma \ref{Le2},
\[
\min\{j\geq1\dvtx\Sigma_j=j-1\}=\min\bigl\{j\geq1\dvtx S^{(+)}_j=-1\bigr\}
=T^{(+)}_1 .
\]
By the ballot theorem [take $f(x^{(\mathrm{c})},x^{(\mathrm{m})})=x^{(\mathrm{c})}+x^{(\mathrm{m})}-1$ in
Lemma \ref{Le6}],
we have
\[
\mathbb{P}(\min\{j\geq1\dvtx\Sigma_j=j-1\}=n\mid\Lambda
_{n-k,k-1})=1/n ,
\]
which yields (i).

The proof of (ii) is similar, observing that
\[
\{|\mathbb{C}_1|=n, M_1=\ell\}=\Lambda_{n-1,\ell}\cap
\bigl\{\min\bigl\{j\geq1\dvtx\Sigma^{(\mathrm{c})}_j=j-1\bigr\}=n\bigr\} .
\]
Finally
(iii) follows by iteration from (ii) and the fact that conditionally on
${A}_1\geq j+1$, the $(j+1)$th allelic cluster $\mathbb{C}_{j+1}$ is independent
of $(\mathbb{C}_k, 1\leq k \leq j)$ and has the same distribution as
the Eve
cluster $\mathbb{C}_1$ (see Lemma \ref{Le3}).
\end{pf}

\subsection{Conditioning on the population size and the number
of alleles}\label{sec32}

In the rest of this section, we will be interested in the relative
sizes of clusters in the allelic partition of the first tree $\mathbb{T}_1$,
ignoring their connections.
We start with a description which is essentially a variation of that in
Theorem \ref{T1}(iii).
Recall that a random uniform cyclic permutation of $\{1,\ldots,k\}$,
say $\sigma$,
is given by $\sigma(i)=U+i$ where~$U$ is uniform on $\{1,\ldots,k\}$
and the addition
is taken modulo $k$.

\begin{corollary}\label{C2} Fix $1\leq k\leq n$ and let $\sigma$ be a
random uniform cyclic permutation of $\{1,\ldots,k\}$ which is
independent of the Galton--Watson process.
Then for every collection of positive integers $n_1,\ldots,n_k$ with
$n_1+\cdots+n_k=n$, we have
\begin{eqnarray*}
&&
\mathbb{P}\bigl(\bigl|\mathbb{C}_{\sigma(1)}\bigr|=n_1,\ldots, \bigl|\mathbb
{C}_{\sigma(k)}\bigr|=n_k
\mid T^{(+)}_1=n, {A}_1=k\bigr)
\\
&&\qquad = \frac{n}{k\pi^{*n}_{n-k,k-1}}
\sum\prod_{i=1}^k
\frac{1}{n_i} \pi^{*n_i}_{n_i-1,\ell_i} ,
\end{eqnarray*}
where in the right-hand side, the sum is taken over the
sequences $\ell_1, \ldots, \ell_k$ in $\mathbb{Z}_+$ such that
$\ell_1+\cdots+\ell_k=k-1$.
\end{corollary}

\begin{pf} A classical application of the ballot theorem shows that
the conditional distribution of $(\mathbb{C}_{\sigma(1)},
\ldots, \mathbb{C}_{\sigma(k)})$ given $T^{(+)}_1=n$ and ${A}_1=k$
is the same as that of $(\mathbb{C}_1, \ldots, \mathbb{C}_k)$
conditioned on
$\sum_{i=1}^{k}|\mathbb{C}_i|=n$ and $\sum_{i=1}^kM_i=k-1$.
Then note that
\[
\sum_{i=1}^{k}|\mathbb{C}_i|=n\quad\mbox{and}\quad\sum_{i=1}^kM_i=k-1
\quad\Longleftrightarrow\quad T^{(\mathrm{c})}_k=n\quad\mbox{and}\quad
\sum_{i=1}^{n}\xi^{(\mathrm{m})}_i=k-1
\]
and an application of the ballot theorem (much in the same way as in
the proof of Theorem \ref{T1}) shows that the probability of that
event equals
\[
\frac{k}{n}\pi^{*n}_{n-k,k-1} .
\]
Theorem \ref{T1}(ii) completes the proof.
\end{pf}

Next, we normalize the size $|\mathbb{C}_i|$
of each cluster by the size $T^{(+)}_1$ of the total population (recall we
focus on the descent of a single ancestor, namely Eve), and write
\[
\Gamma_1\geq\Gamma_2\geq\cdots\geq\Gamma_{{A}_1}
\]
for the sequence which is obtained by ranking the ratios $|\mathbb
{C}_i|/T^{(+)}_1$
in the decreasing order.
So $\Gamma=(\Gamma_1, \ldots, \Gamma_{{A}_1})$ is a proper
partition of
the unit mass, in the sense that it is given by a ranked sequence of
positive real numbers with sum $1$.
The space of mass partitions (possibly with infinitely many strictly
positive terms
and sum less than $1$) is endowed with the supremum distance, which
yields a compact metric space; see Section 2.1 in \cite{RFCP} for details.

Our purpose now is to investigate the asymptotic behavior of
the random mass partition $\Gamma$, under the conditional probability
given the
size $T^{(+)}_{1}=n$ of the tree~$\mathbb{T}_1$ and the number ${A}_1=k$
of alleles,
when $n,k\to\infty$. We shall show that, under appropriate
hypotheses, one
can establish convergence in distribution, where the limit
can be described as follows.
For some fixed parameter $b>0$, consider the sequence
$\mathtt{a}_1>\mathtt{a}_2> \cdots>0$ of the atoms ranked in the decreasing
order of a Poisson point measure on $]0,\infty[$ with intensity
$ba^{-3/2}\,da$.
Roughly speaking, we then get a random proper mass-partition by
conditioning on $\sum_{i=1}^{\infty}\mathtt{a}_i=1$; see, for example,
\cite{PPY} or Proposition 2.4 in \cite{RFCP} for a rigorous definition
of this conditioning by a~singular event.

This family of random mass-partitions has appeared previously in a
remarkable work by Aldous and Pitman \cite{AP1}, more precisely it arose
by logging the Continuum Random Tree according to Poissonian cuts
along its skeleton;
see also \cite{AS,AP2,Ber2,Miermont} for related works. In the
present setting, we may interpret such cuts as mutations which induce
an allelic partition.
As we know from Aldous \cite{Aldous} that the Continuum Random Tree
can be viewed as the
limit when $n\to\infty$ of Galton--Watson trees conditioned to have
total size $n$,
the fact that the preceding random mass-partitions appear again in the
framework of this work should not come as a surprise.

For the sake of simplicity, we shall focus on the case when the number
of clone-children $\xi^{(\mathrm{c})}$ and the number of
mutant-children $\xi^{(\mathrm{m})}$ are independent, although it seems likely that
our argument should also apply to more general situations.
Recall that the expected number of clone-children of a typical
individual is
\mbox{$\mathbb{E}(\xi^{(\mathrm{c})})<1$}. We shall work under the hypothesis
that by a suitable
exponential tilting,
this subcritical random variable can be turned into a critical one with
finite variance.
That is, we shall assume that there exists a real number
$\theta>1$ such that
%
%
\begin{equation}\label{E4}\quad
\mathbb{E}\bigl(\xi^{(\mathrm{c})}\theta^{\xi^{(\mathrm{c})}}\bigr)=\mathbb{E}\bigl(\theta
^{\xi^{(\mathrm{c})}}\bigr) \quad\mbox{and}\quad
\sigma_{\theta}^2:=\mathbb{E}\bigl(\bigl(\xi^{(\mathrm{c})}\bigr)^2 \theta^{\xi^{(\mathrm{c})}}\bigr)/\mathbb
{E}\bigl(\theta^{\xi^{(\mathrm{c})}
}\bigr)-1<\infty.
\end{equation}
It can be readily checked that (\ref{E4}) then specifies $\theta$ uniquely.

\begin{proposition}
Suppose that $\xi^{(\mathrm{c})}$ and $\xi
^{(\mathrm{m})}$ are independent,
that neither distribution is supported by a strict subgroup of $\mathbb{Z}$
and that (\ref{E4}) holds. Fix $b>0$ and let $n,k\to\infty$
according to the regime $k\sim b\sqrt n$. Then the conditional
law of $\Gamma$ given that the size of the total population is $T^{(+)}_1=n$
and the number of alleles ${A}_1=k$ converges weakly on the space of
mass-partitions to the sequence $(\mathtt{a}_1,\mathtt{a}_2, \ldots)$ of the
atoms of a Poisson random measure on $]0,\infty[$ with intensity
\[
\frac{b}{\sqrt{2\pi\sigma_{\theta}^2a^{3}}} \,da ,\qquad
a>0 ,
\]
ranked in the decreasing order and conditioned by $\sum_{i=1}^{\infty
}\mathtt{a}_i=1$.
\end{proposition}

\begin{Remark*}
The special case when $\xi^{(\mathrm{c})}$ and
$\xi^{(\mathrm{m})}$ are
two independent Poisson variables, say with rates $r^{(\mathrm{c})}$ and
$r^{(\mathrm{m})}$ can
also be viewed as an instance
of the situation where mutations affect children independently with probability
$p=r^{(\mathrm{m})}/(r^{(\mathrm{c})}+r^{(\mathrm{m})})$ (cf. the example discussed before Theorem
\ref{T1}).
More precisely
the reproduction law of the standard Galton--Watson process is then
Poisson with rate
$r^{(\mathrm{c})}+r^{(\mathrm{m})}$. This special case has some importance, as it is
well known that
conditioning a Galton--Watson tree with Poisson$(1)$ reproduction law to
have a size $n$
and then assigning to each individual a distinct label in $\{1,\ldots,
n\}$ by uniform sampling without replacements yields the uniform
distribution on the set of rooted trees with $n$ labeled vertices.
\end{Remark*}

\begin{pf} Let $\tilde\mathbb{P}$ denote the probability measure
which is
obtained from $\mathbb{P}$ by exponential tilting, and more precisely,
in such
a way that the variables $\xi^{(\mathrm{c})}_1, \ldots$ are i.i.d. under
$\tilde\mathbb{P}$
with law given by
\[
\tilde\mathbb{P}\bigl(\xi^{(\mathrm{c})}=j\bigr)=\theta^j\mathbb{P}\bigl(\xi^{(\mathrm{c})}=j\bigr)/z_{\theta
} ,\qquad j\in\mathbb{Z}_+ ,
\]
where $z_{\theta}$ is the normalization factor, namely,
\[
z_\theta= \mathbb{E}\bigl(\theta^{\xi^{(\mathrm{c})}}\bigr) .
\]

As in the proof of Corollary \ref{C2}, we see from an application of
the ballot theorem that
the conditional distribution of $(n\Gamma_1,
\ldots, n\Gamma_{{A}_1})$ given $T^{(+)}_1=n$ and ${A}_1=k$
is the same as that obtained from the i.i.d. sequence $|\mathbb{C}_1|,
\ldots,
|\mathbb{C}_k|$ by ranking in the decreasing order and conditioning on
$\sum_{i=1}^{k}|\mathbb{C}_i|=n$ and $\sum_{i=1}^k M_i=k-1$. Observe
that the
latter is
equivalent to conditioning on $\sum_{i=1}^{k}|\mathbb{C}_i|=n$ and
$\sum
_{i=1}^n \xi^{(\mathrm{m})}_i=k-1$.
Further, recall from Lemma \ref{Le3} that $|\mathbb{C}_j|=T^{
(\mathrm{c})}_j-T^{(\mathrm{c})}_{j-1}$ and hence,
on this event, the variables $|\mathbb{C}_1|, \ldots, |\mathbb{C}_k|$
are functions of $\xi^{(\mathrm{c})}_1,\ldots, \xi^{(\mathrm{c})}_n$.
Thus the assumption of independence between
$\xi^{(\mathrm{c})}$ and $\xi^{(\mathrm{m})}$ enables us to ignore the conditioning on
$\sum
_{i=1}^n\xi^{(\mathrm{m})}_i=k-1$.
Finally, it should be clear that the exponential tilting
does not affect such a conditional law, in the sense that the sequence
$|\mathbb{C}_1|, \ldots, |\mathbb{C}_k|$ has the same distribution
under $\mathbb{P}(\cdot
\mid
T^{(+)}_1=n)$
as under $\tilde\mathbb{P}(\cdot\mid T^{(+)}_1=n)$.

We then estimate the distribution of the size of the Eve cluster under
$\tilde\mathbb{P}$, which is given again according to the Dwass formula
\cite
{Dwass} by
\[
\tilde\mathbb{P}(|\mathbb{C}_1|= n_1)=\frac{1}{n_1}\tilde\mathbb
{P}\bigl(\xi^{
(\mathrm{c})}_1+\cdots+ \xi^{(\mathrm{c})}
_{n_1}=n_1-1\bigr)
=\frac{1}{n_1}\tilde\mathbb{P}\bigl(S^{(\mathrm{c})}_{n_1}=-1\bigr) .
\]
Recall that, by assumption, $\xi^{(\mathrm{c})}$ is critical with variance
$\sigma
_{\theta}^2$ under $\tilde\mathbb{P}$,
so an application of Gnedenko's local central limit theorem gives
\[
\tilde\mathbb{P}(|\mathbb{C}_1|= n_1)\sim\frac{1}{\sqrt{2\pi
\sigma_{\theta}^2
n_1^3}} \qquad\mbox{as }n_1\to\infty.
\]

Putting the pieces together, we get that the conditional distribution
of $(n\Gamma_1,
\ldots, n\Gamma_{{A}_1})$ given $T^{(+)}_1=n$ and ${A}_1=k$ is the same
as that
obtained from an i.i.d. sequence\vspace*{1pt} $Y_1, \ldots, Y_k$ by ranking in the
decreasing order and conditioning on
$\sum_{i=1}^{k}Y_i=n$, where
\[
\mathbb{P}(Y_1=n_1)\sim\frac{1}{\sqrt{2\pi\sigma_{\theta}^2
n_1^3}}\qquad
\mbox{as }n_1\to\infty.
\]
An application of Corollary 2.2 in \cite{RFCP} completes the proof of
our claim.
\end{pf}

\section{L\'evy forests with mutations}\label{sec4}

The purpose of this section is to point at an interpretation of a
standard limit theorem involving left-continuous (i.e., downward skip-free)
random walks and L\'evy processes with no negative jumps, in terms of
Galton--Watson and L\'evy forests in the presence of neutral mutations.
We first introduce some notation and hypotheses in this area, referring
to the monograph by Duquesne and Le Gall \cite{DuLG} for details.

For every integer $n\geq1$, let $(\xi^{(\mathrm{c})}(n),\xi^{(\mathrm{m})}(n))$ be a
pair of
integer-valued random variables with
\[
\mathbb{E}\bigl(\xi^{(\mathrm{c})}(n)+\xi^{(\mathrm{m})}(n)\bigr)=1 .
\]
We consider two left-continuous random walks
\[
S^{(+)}(n)=\bigl(S^{(+)}_i(n)\dvtx i\in\mathbb{Z}_+\bigr) \quad\mbox{and}\quad S^{
(\mathrm{c})}(n)=\bigl(S^{(\mathrm{c})}_i(n)\dvtx
i\in\mathbb{Z}_+\bigr),
\]
whose steps are (jointly) distributed as $\xi^{(+)}(n):=\xi^{
(\mathrm{c})}(n)+\xi^{(\mathrm{m})}
(n)-1$ and
$\xi^{(\mathrm{c})}(n)-1$, respectively.
Let also $X=(X_t, t\in\mathbb{R}_+)$ denote
a L\'evy process with no negative jumps
and Laplace exponent $\psi$, namely,
\[
\mathbb{E}(\exp-\lambda X_t)=\exp t\psi(\lambda) \qquad
\mbox{for every }\lambda, t\geq0 .
\]
We further suppose that $X$ does not drift to $+\infty$,
which is equivalent to \mbox{$\psi'(0+)\geq0$}, and that
\[
\int_1^{\infty}\frac{d\lambda}{\psi(\lambda)}<\infty.
\]

We also need to introduce a different procedure for encoding forests by
paths, which
is more convenient to work with when discussing continuous limits of
discrete structures. For each $n\geq1$, we write
$H(n)=(H_i(n), i\in\mathbb{N})$ for the (discrete) \textit{height
function} of the
Galton--Watson forest
$(\mathbb{T}_{\ell}, \ell\in\mathbb{N})$. That is, for $i\geq0$,
$H_i(n)$ denotes the generation of the $(i+1)$th individual found by the
usual depth-first search (i.e., mutations are discarded)
on the Galton--Watson forest. In the continuous setting, trees and
forests can be defined for a fairly general class of L\'evy processes
with no negative jumps, and in turn are encoded by (continuous) height
functions; cf. Chapter 1 in \cite{DuLG} for precise definitions and
further references.

The key hypothesis in this setting is the existence of a nondecreasing
sequence of positive integers $(\gamma_n, n\in\mathbb{N})$
converging to $\infty$ and such that
%
%
\begin{equation}\label{E6}
\lim_{n\to\infty} n^{-1} S^{(+)}_{n\gamma_n}(n)=X_1\qquad \mbox{in law};
\end{equation}
we also assume that the technical condition \textup{(2.27)} in \cite{DuLG}
is fulfilled. Then
the rescaled height function
\[
\bigl(\gamma_n^{-1}H_{[tn\gamma_n]}(n)\dvtx t\geq0 \bigr)
\]
converges in distribution, in the sense of weak convergence on Skorohod
space $\mathbb{D}(\mathbb{R}_+,\mathbb{R}_+)$ as $n\to\infty$
toward the height process
$ ({\mathcal H}_t\dvtx t\geq0 )$
which is constructed from the L\'evy process $X=(X_t, t\geq0)$; see
Theorem 2.3.1 in \cite{DuLG}.

Similarly, we write $H^{(\mathrm{c})}(n)=(H^{(\mathrm{c})}_i(n), i\in\mathbb
{N})$ for
the height function
of the Galton--Watson forest $(\mathbb{C}_{j}, j\in\mathbb{N})$,
where each allelic
cluster $\mathbb{C}_j$ is endowed with the genealogical tree structure induced
by the population model (see Remark, item~1 in Section \ref{sec23}).

\begin{proposition} \label{P2} Suppose that the preceding assumptions
hold, and also that
%
%
\begin{equation}\label{E7}
\lim_{n\to\infty} \gamma_n \mathbb{E}\bigl(\xi^{(\mathrm{m})}(n)\bigr)=d
\quad\mbox{and}\quad \lim_{n\to\infty} n^{-1}\gamma_n\operatorname{Var}\bigl(\xi^{(\mathrm{m})}(n)\bigr)=0
\end{equation}
for some $d\geq0$. Then
the rescaled height function
\[
\bigl(\gamma_n^{-1}H^{(\mathrm{c})}_{[tn \gamma_n]}(n) \dvtx t\geq0 \bigr)
\]
converges in distribution, in the sense of weak convergence on Skorohod
space $\mathbb{D}(\mathbb{R}_+,\mathbb{R}_+)$ as $n\to\infty$
toward the height process
\[
\bigl( {\mathcal H}_t^{(d)}\dvtx t\geq0 \bigr) ,
\]
which is constructed from the L\'evy process
$X^{(d)}=(X^{(d)}_t:=X_t-dt, t\geq0)$.
\end{proposition}

\begin{Remark*} More recently, Duquesne and Le Gall \cite{DuLG2}
(see also the survey~\cite{LG}) have developed the framework when L\'
evy trees are viewed as random variables with values in the space of
real trees, endowed with the Gromov--Hausdorff distance. Proposition
\ref
{P2} can also be restated in this setting.
\end{Remark*}

\begin{pf*}{Proof of Proposition \protect\ref{P2}}
The assumption (\ref{E6}) ensures the convergence in distribution
\[
\bigl(n^{-1}S^{(+)}_{[tn \gamma_n]}(n)\dvtx t\geq0 \bigr)
\Longrightarrow(X_t\dvtx t\geq0 ) ,
\]
see Theorem 2.1.1 in \cite{DuLG}
and (2.3) there. On the other hand, by a routine argument
based on martingales, the assumption (\ref{E7}) entails that
\[
\lim_{n\to\infty} n^{-1}\bigl(S^{(+)}_{[tn\gamma_n]}(n)-S^{
(\mathrm{c})}_{[tn\gamma_n]}(n)\bigr)=dt,
\]
uniformly for $t$ in compact intervals, in $L^2(\mathbb{P})$.
The convergence in distribution
\[
\bigl(n^{-1}S^{(\mathrm{c})}_{[tn\gamma_n]}(n) \dvtx t\geq0 \bigr)
\Longrightarrow(X_t-dt\dvtx t\geq0 )
\]
follows. Recall that depth-first search with mutations on the initial
forest yields the usual depth-first search for the forest of allelic
clusters (cf. Remark, item 1 in Section~\ref{sec23}).
We can then complete the proof as in Theorem 2.3.1 in \cite{DuLG}.
\end{pf*}

We now conclude this work by discussing a natural example.
Specifically, we suppose that the distribution of
\[
\xi^{(\mathrm{c})}(n)+\xi^{(\mathrm{m})}(n)=\xi(n):= \xi
\]
is
the same for all $n$.
For the sake of simplicity, we assume also that $\mathbb{E}(\xi)=1$
and $\operatorname{Var}(\xi)=1$. We may then take $\gamma_n=n$, so by
the central limit
theorem, (\ref{E6}) holds and the L\'evy process $X$ is a standard
Brownian motion. We fix an arbitrary $d>0$ and consider the independent
pruning model where for each integer $n>d$,
conditionally on the total number of children $\xi^{(+)}(n):=\xi^{
(\mathrm{c})}(n)+\xi^{(\mathrm{m})}(n)=k$,
the number $\xi^{(\mathrm{m})}(n)$ of mutant-children of a typical individual has
the binomial distribution $B(k, d/n)$.
In other words, in the $n$th population model, mutations affect each
child with probability $d/n$, independently of the other children.
Then~(\ref{E7}) clearly holds.
Roughly speaking, Theorem 2.3.1 of \cite{DuLG} implies in this setting that
the initial Galton--Watson forest associated with the $n$th population
model, converges
in law after a suitable renormalization to the Brownian forest, whereas
Proposition~\ref{P2} of the present work shows that the allelic forest
renormalized in the same way, converges in law to the forest generated
by a Brownian motion with drift~$-d$.\looseness=1

This provides an explanation to the rather intriguing relation which
identifies two seemingly different fragmentation processes: the
fragmentation process constructed by Aldous and Pitman \cite{AP1} by
logging the Continuum Random Tree according to a Poisson point process
on its skeleton,
and the fragmentation process constructed in \cite{Ber2} by splitting
the unit interval at instants when
a Brownian excursion with negative drift reaches a new infimum.
It is interesting to mention that Schweinsberg \cite{Schw} already
pointed at several applications of the (continuous) ballot theorem in
this framework.
More generally, the transformation $X\to X^{(d)}$ of L\'evy processes
with no negative jumps also appeared in an article by Miermont \cite
{Miermont} on certain eternal additive coalescents, whereas
Aldous and Pitman \cite{AP2} showed that the latter arise
asymptotically from independent pruning of certain sequences of
birthday trees.
Finally, we also refer \cite{ADV} for another interesting recent work
on pruning L\'evy random trees.

\section*{Acknowledgment}
I would like to thank two anonymous
referees for their careful check of this work.

%

%
\printaddresses


\begin{thebibliography}{40}

\bibitem{AD}
%
\begin{bmisc}[unstr]
\bauthor{\bsnm{Abraham},~\bfnm{R.}\binits{R.}} \AND
\bauthor{\bsnm{Delmas},~\bfnm{J.-F.}\binits{J.-F.}}
(\byear{2008}).
\bhowpublished{Williams' decomposition of the L\'evy continuous
random tree and simultaneous extinction probability for populations with
neutral mutations. Preprint}.
Available at \url{http://hal.archives-ouvertes.fr/}.
\end{bmisc}
%
\endbibitem

\bibitem{ADV}
%
\begin{bmisc}[unstr]
\bauthor{\bsnm{Abraham},~\bfnm{R.}\binits{R.}},
\bauthor{\bsnm{Delmas},~\bfnm{J.-F.}\binits{J.-F.}} \AND
\bauthor{\bsnm{Voisin},~\bfnm{G.}\binits{G.}}
(\byear{2009}).
Pruning a L\'{e}vy continuum random tree.
\textit{Stochastic Process. Appl.} To appear.
\end{bmisc}
%
\endbibitem

\bibitem{AS}
%
\begin{barticle}[vtex]
\bauthor{\bsnm{Abraham},~\bfnm{Romain}\binits{R.}} \AND
\bauthor{\bsnm{Serlet},~\bfnm{Laurent}\binits{L.}}
(\byear{2002}).
\btitle{Poisson snake and fragmentation}.
\bjournal{Electron. J. Probab.}
\bvolume{7}
\bpages{15 (electronic)}.
\bmrnumber{MR1943890}
\end{barticle}
%
\endbibitem

\bibitem{Aldous}
%
\begin{barticle}[msn]
\bauthor{\bsnm{Aldous},~\bfnm{David}\binits{D.}}
(\byear{1993}).
\btitle{The continuum random tree. {III}}.
\bjournal{Ann. Probab.}
\bvolume{21}
\bpages{248--289}.
\bmrnumber{MR1207226}
\end{barticle}
%
\endbibitem

\bibitem{AP0}
%
\begin{barticle}[vtex]
\bauthor{\bsnm{Aldous},~\bfnm{David}\binits{D.}} \AND
\bauthor{\bsnm{Pitman},~\bfnm{Jim}\binits{J.}}
(\byear{1998}).
\btitle{Tree-valued {M}arkov chains derived from {G}alton--{W}atson processes}.
\bjournal{Ann. Inst. H. Poincar\'e Probab. Statist.}
\bvolume{34}
\bpages{637--686}.
\bmrnumber{MR1641670}
\end{barticle}
%
\endbibitem

\bibitem{AP1}
%
\begin{barticle}[msn]
\bauthor{\bsnm{Aldous},~\bfnm{David}\binits{D.}} \AND
\bauthor{\bsnm{Pitman},~\bfnm{Jim}\binits{J.}}
(\byear{1998}).
\btitle{The standard additive coalescent}.
\bjournal{Ann. Probab.}
\bvolume{26}
\bpages{1703--1726}.
\bmrnumber{MR1675063}
\end{barticle}
%
\endbibitem

\bibitem{AP2}
%
\begin{barticle}[msn]
\bauthor{\bsnm{Aldous},~\bfnm{David}\binits{D.}} \AND
\bauthor{\bsnm{Pitman},~\bfnm{Jim}\binits{J.}}
(\byear{2000}).
\btitle{Inhomogeneous continuum random trees and the entrance boundary
of the
additive coalescent}.
\bjournal{Probab. Theory Related Fields}
\bvolume{118}
\bpages{455--482}.
\bmrnumber{MR1808372}
\end{barticle}
%
\endbibitem

\bibitem{AN}
%
\begin{bbook}[vtex]
\bauthor{\bsnm{Athreya},~\bfnm{Krishna~B.}\binits{K.~B.}} \AND
\bauthor{\bsnm{Ney},~\bfnm{Peter~E.}\binits{P.~E.}}
(\byear{1972}).
\btitle{Branching Processes}.
\bpublisher{Springer}, \baddress{New York}.
\bmrnumber{MR0373040}
\end{bbook}
%
\endbibitem

\bibitem{BG}
%
\begin{barticle}[vtex]
\bauthor{\bsnm{Basdevant},~\bfnm{Anne-Laure}\binits{A.-L.}} \AND
\bauthor{\bsnm{Goldschmidt},~\bfnm{Christina}\binits{C.}}
(\byear{2008}).
\btitle{Asymptotics of the allele frequency spectrum associated with the
{B}olthausen--{S}znitman coalescent}.
\bjournal{Electron. J. Probab.}
\bvolume{13}
\bpages{486--512}.
\bmrnumber{MR2386740}
\end{barticle}
%
\endbibitem

\bibitem{BBS}
%
\begin{barticle}[msn]
\bauthor{\bsnm{Berestycki},~\bfnm{Julien}\binits{J.}},
\bauthor{\bsnm{Berestycki},~\bfnm{Nathana{\"e}l}\binits{N.}} \AND
\bauthor{\bsnm{Schweinsberg},~\bfnm{Jason}\binits{J.}}
(\byear{2007}).
\btitle{Beta-coalescents and continuous stable random trees}.
\bjournal{Ann. Probab.}
\bvolume{35}
\bpages{1835--1887}.
\bmrnumber{MR2349577}
\end{barticle}
%
\endbibitem

\bibitem{Ber1}
%
\begin{barticle}[msn]
\bauthor{\bsnm{Bertoin},~\bfnm{Jean}\binits{J.}}
(\byear{1999}).
\btitle{Renewal theory for embedded regenerative sets}.
\bjournal{Ann. Probab.}
\bvolume{27}
\bpages{1523--1535}.
\bmrnumber{MR1733158}
\end{barticle}
%
\endbibitem

\bibitem{Ber2}
%
\begin{barticle}[msn]
\bauthor{\bsnm{Bertoin},~\bfnm{Jean}\binits{J.}}
(\byear{2000}).
\btitle{A fragmentation process connected to {B}rownian motion}.
\bjournal{Probab. Theory Related Fields}
\bvolume{117}
\bpages{289--301}.
\bmrnumber{MR1771665}
\end{barticle}
%
\endbibitem

\bibitem{RFCP}
%
\begin{bbook}[msn]
\bauthor{\bsnm{Bertoin},~\bfnm{Jean}\binits{J.}}
(\byear{2006}).
\btitle{Random Fragmentation and Coagulation Processes}.
\bseries{Cambridge Studies in Advanced Mathematics}
\bvolume{102}.
\bpublisher{Cambridge Univ. Press}, \baddress{Cambridge}.
\bmrnumber{MR2253162}
\end{bbook}
%
\endbibitem

\bibitem{CD}
%
\begin{barticle}[vtex]
\bauthor{\bsnm{Crump},~\bfnm{Kenny~S.}\binits{K.~S.}} \AND
\bauthor{\bsnm{Gillespie},~\bfnm{John~H.}\binits{J.~H.}}
(\byear{1976}).
\btitle{The dispersion of a neutral allele considered as a branching process}.
\bjournal{J. Appl. Probab.}
\bvolume{13}
\bpages{208--218}.
\bmrnumber{MR0408877}
\end{barticle}
%
\endbibitem

\bibitem{DDSJ}
%
\begin{barticle}[msn]
\bauthor{\bsnm{Delmas},~\bfnm{Jean-Fran\c{c}ois}\binits{J.-F.}},
\bauthor{\bsnm{Dhersin},~\bfnm{Jean-St{\'e}phane}\binits{J.-S.}}
\AND
\bauthor{\bsnm{Siri-Jegousse},~\bfnm{Arno}\binits{A.}}
(\byear{2008}).
\btitle{Asymptotic results on the length of coalescent trees}.
\bjournal{Ann. Appl. Probab.}
\bvolume{18}
\bpages{997--1025}.
\bmrnumber{MR2418236}
\end{barticle}
%
\endbibitem

\bibitem{DGP}
%
\begin{barticle}[msn]
\bauthor{\bsnm{Dong},~\bfnm{Rui}\binits{R.}},
\bauthor{\bsnm{Gnedin},~\bfnm{Alexander}\binits{A.}} \AND
\bauthor{\bsnm{Pitman},~\bfnm{Jim}\binits{J.}}
(\byear{2007}).
\btitle{Exchangeable partitions derived from {M}arkovian coalescents}.
\bjournal{Ann. Appl. Probab.}
\bvolume{17}
\bpages{1172--1201}.
\bmrnumber{MR2344303}
\end{barticle}
%
\endbibitem

\bibitem{DuLG}
%
\begin{barticle}[vtex]
\bauthor{\bsnm{Duquesne},~\bfnm{Thomas}\binits{T.}} \AND
\bauthor{\bsnm{Le~Gall},~\bfnm{Jean-Fran\c{c}ois}\binits{J.-F.}}
(\byear{2002}).
\btitle{Random trees, {L}\'evy processes and spatial branching processes}.
\bjournal{Ast\'erisque}
\bvolume{281}
\bpages{vi--147}.
\bmrnumber{MR1954248}
\end{barticle}
%
\endbibitem

\bibitem{DuLG2}
%
\begin{barticle}[vtex]
\bauthor{\bsnm{Duquesne},~\bfnm{Thomas}\binits{T.}} \AND
\bauthor{\bsnm{Le~Gall},~\bfnm{Jean-Fran\c{c}ois}\binits{J.-F.}}
(\byear{2005}).
\btitle{Probabilistic and fractal aspects of {L}\'evy trees}.
\bjournal{Probab. Theory Related Fields}
\bvolume{131}
\bpages{553--603}.
\bmrnumber{MR2147221}
\end{barticle}
%
\endbibitem

\bibitem{Dwass}
%
\begin{barticle}[msn]
\bauthor{\bsnm{Dwass},~\bfnm{Meyer}\binits{M.}}
(\byear{1969}).
\btitle{The total progeny in a branching process and a related random walk.}
\bjournal{J.~Appl. Probab.}
\bvolume{6}
\bpages{682--686}.
\bmrnumber{MR0253433}
\end{barticle}
%
\endbibitem

\bibitem{Ewens}
%
\begin{barticle}[vtex]
\bauthor{\bsnm{Ewens},~\bfnm{W.~J.}\binits{W.~J.}}
(\byear{1972}).
\btitle{The sampling theory of selectively neutral alleles}.
\bjournal{Theoret. Population Biology}
\bvolume{3}
\bpages{87--112}.
\bmrnumber{MR0325177}
\end{barticle}
%
\endbibitem

\bibitem{Fel1}
%
\begin{bbook}[vtex]
\bauthor{\bsnm{Feller},~\bfnm{William}\binits{W.}}
(\byear{1968}).
\btitle{An Introduction to Probability Theory and Its Applications}
\bvolume{I}.
\bpublisher{Wiley}, \baddress{New York}.
\bmrnumber{MR0228020}
\end{bbook}
%
\endbibitem

\bibitem{Fel2}
%
\begin{bbook}[vtex]
\bauthor{\bsnm{Feller},~\bfnm{William}\binits{W.}}
(\byear{1971}).
\btitle{An Introduction to Probability Theory and Its Applications}
\bvolume{II}.
\bpublisher{Wiley}, \baddress{New York}.
\bmrnumber{MR0270403}
\end{bbook}
%
\endbibitem

\bibitem{GP}
%
\begin{barticle}[msn]
\bauthor{\bsnm{Griffiths},~\bfnm{R.~C.}\binits{R.~C.}} \AND
\bauthor{\bsnm{Pakes},~\bfnm{Anthony~G.}\binits{A.~G.}}
(\byear{1988}).
\btitle{An infinite-alleles version of the simple branching process}.
\bjournal{Adv. in Appl. Probab.}
\bvolume{20}
\bpages{489--524}.
\bmrnumber{MR955502}
\end{barticle}
%
\endbibitem

\bibitem{Harris1}
%
\begin{barticle}[msn]
\bauthor{\bsnm{Harris},~\bfnm{T.~E.}\binits{T.~E.}}
(\byear{1952}).
\btitle{First passage and recurrence distributions}.
\bjournal{Trans. Amer. Math. Soc.}
\bvolume{73}
\bpages{471--486}.
\bmrnumber{MR0052057}
\end{barticle}
%
\endbibitem

\bibitem{Harris2}
%
\begin{bbook}[vtex]
\bauthor{\bsnm{Harris},~\bfnm{Theodore~E.}\binits{T.~E.}}
(\byear{1963}).
\btitle{The Theory of Branching Processes}.
\bpublisher{Springer}, \baddress{Berlin}.
\bmrnumber{MR0163361}
\end{bbook}
%
\endbibitem

\bibitem{KimmAxel}
%
\begin{bbook}[msn]
\bauthor{\bsnm{Kimmel},~\bfnm{Marek}\binits{M.}} \AND
\bauthor{\bsnm{Axelrod},~\bfnm{David~E.}\binits{D.~E.}}
(\byear{2002}).
\btitle{Branching Processes in Biology}.
\bseries{Interdisciplinary Applied Mathematics}
\bvolume{19}.
\bpublisher{Springer}, \baddress{New York}.
\bmrnumber{MR1903571}
\end{bbook}
%
\endbibitem

\bibitem{Kingman}
%
\begin{bbook}[vtex]
\bauthor{\bsnm{Kingman},~\bfnm{J.~F.~C.}\binits{J.~F.~C.}}
(\byear{1980}).
\btitle{Mathematics of Genetic Diversity}.
\bseries{CBMS-NSF Regional Conference Series in Applied Mathematics}
\bvolume{34}.
\bpublisher{SIAM},
\baddress{Philadelphia, PA.}
\bmrnumber{MR591166}
\end{bbook}
%
\endbibitem



\bibitem{Lambert}
%
\begin{bmisc}[unstr]
\bauthor{\bsnm{Lambert},~\bfnm{A.}\binits{A.}}
(\byear{2008}).
\bhowpublished{Spine decompositions and allelic partitions of
splitting trees. In preparation}.
\end{bmisc}
%
\endbibitem

\bibitem{LG}
%
\begin{barticle}[msn]
\bauthor{\bsnm{Le~Gall},~\bfnm{Jean-Fran\c{c}ois}\binits{J.-F.}}
(\byear{2005}).
\btitle{Random trees and applications}.
\bjournal{Probab. Surv.}
\bvolume{2}
\bpages{245--311 (electronic)}.
\bmrnumber{MR2203728}
\end{barticle}
%
\endbibitem

\bibitem{LSS}
%
\begin{barticle}[msn]
\bauthor{\bsnm{Liggett},~\bfnm{Thomas~M.}\binits{T.~M.}},
\bauthor{\bsnm{Schinazi},~\bfnm{Rinaldo~B.}\binits{R.~B.}} \AND
\bauthor{\bsnm{Schweinsberg},~\bfnm{Jason}\binits{J.}}
(\byear{2008}).
\btitle{A contact process with mutations on a tree}.
\bjournal{Stochastic Process. Appl.}
\bvolume{118}
\bpages{319--332}.
\bmrnumber{MR2389047}
\end{barticle}
%
\endbibitem

\bibitem{LyonsPeres}
%
\begin{bmisc}[unstr]
\bauthor{\bsnm{Lyons},~\bfnm{R.}\binits{R.}} \AND
\bauthor{\bsnm{Peres},~\bfnm{Y.}\binits{Y.}}
(\byear{2008}).
\bhowpublished{Probability on trees and networks}.
Available at
\href{http://php.indiana.edu/~rdlyons/prbtree/book.pdf}{http://php.indiana.edu/\texttildelow rdlyons/prbtree/book.pdf}.
\end{bmisc}
%
\endbibitem

\bibitem{Miermont}
%
\begin{barticle}[vtex]
\bauthor{\bsnm{Miermont},~\bfnm{Gr{\'e}gory}\binits{G.}}
(\byear{2001}).
\btitle{Ordered additive coalescent and fragmentations associated to {L}\'{e}vy
processes with no positive jumps}.
\bjournal{Electron. J. Probab.}
\bvolume{6}
\bpages{33 (electronic)}.
\bmrnumber{MR1844511}
\end{barticle}
%
\endbibitem

\bibitem{Moehle}
%
\begin{barticle}[msn]
\bauthor{\bsnm{M{\"o}hle},~\bfnm{M.}\binits{M.}}
(\byear{2006}).
\btitle{On sampling distributions for coalescent processes with simultaneous
multiple collisions}.
\bjournal{Bernoulli}
\bvolume{12}
\bpages{35--53}.
\bmrnumber{MR2202319}
\end{barticle}
%
\endbibitem

\bibitem{Nerman}
%
\begin{binproceedings}[vtex]
\bauthor{\bsnm{Nerman},~\bfnm{O.}\binits{O.}}
(\byear{1987}).
\btitle{Branching processes and neutral mutations}.
In \bbooktitle{Proceedings of the First {W}orld {C}ongress of the {B}ernoulli
{S}ociety, {V}ol. 2 ({T}ashkent, 1986)}
\bpages{683--692}.
\bpublisher{VNU Sci. Press}, \baddress{Utrecht}.
\bmrnumber{MR1092502}
\end{binproceedings}
%
\endbibitem

\bibitem{PPY}
%
\begin{barticle}[msn]
\bauthor{\bsnm{Perman},~\bfnm{Mihael}\binits{M.}},
\bauthor{\bsnm{Pitman},~\bfnm{Jim}\binits{J.}} \AND
\bauthor{\bsnm{Yor},~\bfnm{Marc}\binits{M.}}
(\byear{1992}).
\btitle{Size-biased sampling of {P}oisson point processes and excursions}.
\bjournal{Probab. Theory Related Fields}
\bvolume{92}
\bpages{21--39}.
\bmrnumber{MR1156448}
\end{barticle}
%
\endbibitem

\bibitem{PiSF}
%
\begin{bbook}[vtex]
\bauthor{\bsnm{Pitman},~\bfnm{J.}\binits{J.}}
(\byear{2006}).
\btitle{Combinatorial Stochastic Processes}.
\bseries{Lecture Notes in Math.}
\bvolume{1875}.
\bpublisher{Springer}, \baddress{Berlin}.
\bmrnumber{MR2245368}
\end{bbook}
%
\endbibitem

\bibitem{SS}
%
\begin{barticle}[vtex]
\bauthor{\bsnm{Schinazi},~\bfnm{R.~B.}\binits{R.~B.}} \AND
\bauthor{\bsnm{Schweinsberg},~\bfnm{J.}\binits{J.}}
(\byear{2008}).
\btitle{Spatial and nonspatial stochastic models for immune response}.
\bjournal{Markov Process. Related Fields}
\bvolume{14}
\bpages{255--276}.
\bmrnumber{MR2437531}
\end{barticle}
%
\endbibitem

\bibitem{Schw}
%
\begin{barticle}[msn]
\bauthor{\bsnm{Schweinsberg},~\bfnm{Jason}\binits{J.}}
(\byear{2001}).
\btitle{Applications of the continuous-time ballot theorem to
{B}rownian motion
and related processes}.
\bjournal{Stochastic Process. Appl.}
\bvolume{95}
\bpages{151--176}.
\bmrnumber{MR1847096}
\end{barticle}
%
\endbibitem


\bibitem{Taib}
%
\begin{bbook}[msn]
\bauthor{\bsnm{Ta{\"{\i}}b},~\bfnm{Ziad}\binits{Z.}}
(\byear{1992}).
\btitle{Branching Processes and Neutral Evolution}.
\bseries{Lecture Notes in Biomathematics}
\bvolume{93}.
\bpublisher{Springer}, \baddress{Berlin}.
\bmrnumber{MR1176317}
\end{bbook}
%
\endbibitem

\bibitem{Takacs}
%
\begin{bbook}[msn]
\bauthor{\bsnm{Tak{\'a}cs},~\bfnm{Lajos}\binits{L.}}
(\byear{1967}).
\btitle{Combinatorial Methods in the Theory of Stochastic Processes}.
\bpublisher{Wiley}, \baddress{New York}.
\bmrnumber{MR0217858}
\end{bbook}
%
\endbibitem

\end{thebibliography}
\end{document}